\documentclass[12pt]{amsart}
\title{Hurwitz Theory and the Double Ramification Cycle}
\author{Renzo Cavalieri}
\address{Renzo Cavalieri, Colorado State University, Department of Mathematics, Weber Building, Fort Collins, CO 80523-1874, USA}
\email{renzo@math.colostate.edu}

\usepackage {color,graphicx,psfrag,enumitem, verbatim,amssymb,color,xypic}
\usepackage{texdraw}

\newcommand{\proj}{\mathbb{P}}
\newcommand{\R}{\mathbb{R}}
\newcommand{\N}{\mathbb{N}}
\newcommand{\Z}{\mathbb{Z}}
\newcommand{\aut}{\mathrm{Aut}}

\newcommand{\M}{\mathcal{M}}
\newcommand{\so}{\mathcal{O}}
\newcommand{\Q}{\mathbb{Q}}
\newcommand{\C}{\mathbb{C}}
\newcommand{\sE}{\mathcal{E}}
\newcommand{\sO}{\mathcal{O}}
\newcommand{\fC}{\mathcal{C}}
\newcommand{\sM}{\mathcal{M}}
\newcommand{\sU}{\mathcal{U}}
\newcommand{\sD}{\mathcal{D}}
\newcommand{\bZ}{\mathbb{Z}}

\newcommand{\cham}{\mathfrak{c}} 
\newcommand{\A}{\mathcal{A}}  
\newcommand{\ch}{\mbox{ch}} 

\newcommand{\minus}{\ \mbox{\footnotesize{$\setminus$}}\ }

\newcommand{\cms}{\overline{\sM}} 

\newcommand{\relrub}[2]{\cms^\sim_{#1}(\PP^1;#2)}

\newcommand{\bx}{\textbf{x}}

\newcommand{\beq}{\begin{equation}}
\newcommand{\eeq}{\end{equation}}
\newcommand{\PP}{\mathbb{P}}    
\providecommand{\stab}{{\rm stab}}

\newcommand{\br}{{\rm br}}

\newcommand{\pt}{{\rm pt}}

\newtheorem{dummy}{}[section]
\newtheorem{lemma}[dummy]{Lemma}

\newtheorem{theorem}[dummy]{Theorem}
\newtheorem{corollary}[dummy]{Corollary}
\newtheorem*{conjecture}{Conjecture}
\newtheorem*{question}{Question}

\theoremstyle{definition}
\newtheorem{definition}[dummy]{Definition}
\newtheorem{example}[dummy]{Example}
\newtheorem*{fact}{Fact}
\newtheorem{remark}[dummy]{Remark}
\newtheorem{exercise}{Exercise}
\newcommand{\LM}{\overline{M}_{0,2+r}(1,1,\varepsilon, \ldots, \varepsilon)}
\newcommand{\mon}{\overline{M}_{0,n}}
\newcommand{\Htk}{\tilde\bbH_k(\mathbf{x})}

\newcommand{\bbH}{\mathbb{H}}
\newcommand{\fc}{\mathfrak{c}}

\newcommand{\cT}{\mathcal{T}}
\newcommand{\bC}{\mathbb{C}}
\newcommand{\cZ}{\mathcal{Z}}

\def\cH{{\mathcal{H}}}

\def\cM{{\mathcal{M}}}

\def\cT{{\cal T}}

\def\cZ{{\mathcal{Z}}}

\newcommand{\Mbar}{\overline{\cM}}
\newcommand{\Hbar}{\overline{\cH}}

1
\begin{document}
\maketitle
\begin{abstract}
This survey grew out of notes  accompanying a cycle of lectures at the workshop {\it Modern Trends in Gromov-Witten Theory}, in Hannover. The lectures are devoted to  interactions between Hurwitz theory and Gromov-Witten theory, with a particular eye to the contributions made to the understanding of the Double Ramification Cycle, a cycle in the moduli space of curves that compactifies the double Hurwitz locus.
We explore the algebro-combinatorial properties of single and double Hurwitz numbers, and the connections with  intersection theoretic problems on appropriate moduli spaces. 
We survey several results by many groups of people on the subject, but, perhaps more importantly, collect  a  number of conjectures  and problems which are still open.

\end{abstract}

\tableofcontents

 \section*{Introduction}
 This article surveys a series of mathematical works and ideas centered on the interactions between Hurwitz theory, Gromov-Witten theory and tautological intersection theory on the moduli space of curves. After a skeletal introduction to Hurwitz theory, we focus on geometric and combinatorial properties of specific families of Hurwitz numbers, and the corresponding geometric loci, representing interesting tautological cycles in the moduli space of curves.
 
We begin with the ELSV formula \cite{elsv:hnaiomsoc}, which expresses single Hurwitz numbers as  linear Hodge integrals, intersection numbers of certain tautological classes on the moduli space  of curves. This formula has remarkably been fertile in ``both directions": on the one hand, Hurwitz numbers are readily computable, and hence the ELSV formula has aided in the computations of Hodge integrals. On the other, the ELSV formula exhibits single Hurwitz numbers (up to some combinatorial prefactor) as polynomials in the ramification data with coefficients Hodge monomials, thus shedding light on a structure for single Hurwitz numbers observed and conjectured by Goulden and Jackson \cite{gj:tf}.

Next we focus on double Hurwitz numbers, where a similar structure was observed by Goulden, Jackson and Vakil: families of double Hurwitz numbers are piecewise polynomial in the ramification data. In this case the piecewise polynomial structure and  wall crossing formulae  were proven combinatorially in \cite{gjv:ttgodhn,ssv:g0,cjm:thn,cjm:wcfdhn}, giving rise to a number of geometric questions and conjectures. One natural objective, originally in \cite{gjv:ttgodhn}, is the search for an ``ELSV-type" formula expressing double Hurwitz numbers as intersection numbers over a family of compactifications of the Picard stack. We present here a generalized version of Goulden, Jackson and Vakil's original conjecture; our generalization attempts to explain the piecewise polynomial structure as a variation of stability on a family of moduli spaces. While finding an ELSV formula on a family of Picard stacks is still an open question, in recent work with Marcus \cite{cm:dhn} we exhibit an intersection theoretic formula for double Hurwitz numbers which explains the wall-crossings as a consequence of chamber dependent corrections to the $\psi$ classes on $\overline{M}_{g,n}$.

Developing such a formula requires working with the moduli spaces of covers with discrete data specified by the double Hurwitz number, and a suitable compactfication thereof which allows us to carry on intersection theory. In fact, the search for a cycle which meaningfully compactifies  the locus of curves  supporting a map with specified ramification data is known as {\it Eliashberg's question}; in the last ten or so years several groups of mathematicians have been contributing to our, as of yet incomplete, understanding of such a cycle, which is now known as the double ramification cycle (\cite{H11, bssz:psi,GZ1,GZ2}). The most exciting development in this story is a recent conjecture of Aaron Pixton, expressing explicitly the class of the double ramification cycle as a combination of standard tautological classes. We conclude this survey by discussing Pixton's conjecture.

Finally,  Section \ref{trop} gives a brief survey of the role that tropical geometry is playing in this story. Tropical covers have from the very beginning been an extremely efficient way of organizing the count of Hurwitz covers (see also Johnson's survey \cite{Paul:hnrgt}). Recent work of Abramovich, Caporaso and Payne  \cite{ACP} has put the connection between tropical and classical moduli spaces of curves on firm conceptual ground, exhibiting the tropical moduli space as the Berkovich skeleton  of the analytification of the moduli space of curves, associated to the natural toroidal structure of $\overline{M}_{g,n}$. In recent work with Markwig and Ranganathan \cite{CMR14,CMR14b}, we extend this point of view to moduli spaces of admissible covers and relative stable maps to $\PP^1$, which allows us exhibit the classical/tropical correspondence as a geometric fact, rather than a combinatorial one.

We are intending  these notes to serve as a useful roadmap for a graduate student beginning his exploration of the field, or for a researcher who has been interested in some aspects of this story and is seeking a more complete overview.

\subsection*{Acknowledgements}
I am grateful to the many people that I have worked with and talked to about these topics, including my collaborators Aaron Bertram, Paul Johnson, Steffen Marcus, Hannah Markwig, Dhruv Ranganathan and Jonathan Wise; special thanks to Arend Bayer, Brian Osserman and Dimitri Zvonkine: many discussions with them have deeply contributed to shape my point of view on the subject.
I acknowledge with gratitude the NSF's support through grant DMS-1101549, RTG grant 1159964; thanks to the organizers of the summer school {\it Modern Trends in Gromov-Witten Theory} in Hannover for the opportunity to give this cycle of lectures. A great thanks to the anonymous referee who gave a long list of excellent comments, and caused me to improve the quality of the exposition. Thanks to my graduate student Vance Blankers for helping me catch a few typos, especially some embarrassing ones in formulas.

\section{From Hurwitz to ELSV}
\begin{figure}[b]
	\centering
		\includegraphics[width=0.75\textwidth]{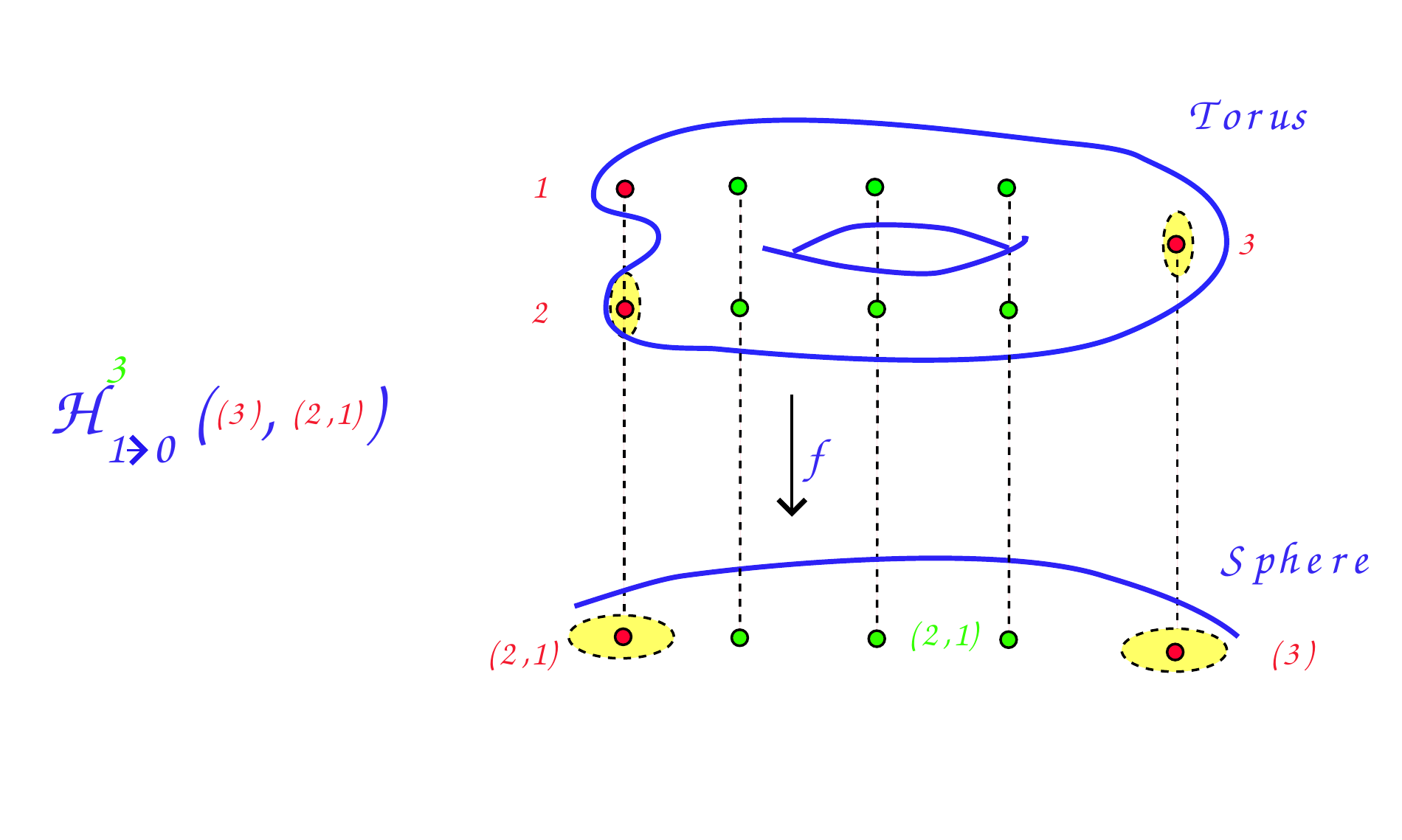}
	\caption{The covers contributing to a given Hurwitz Number.}
	\label{Hurw}
\end{figure}

\begin{definition}[Geometry]\label{hngeom}
Let $ (Y,p_1,\ldots,p_r,q_1,\ldots,q_s) $ be an $(r+s)$-marked smooth Riemann Surface of genus $g$.  Let $
\underline{\eta}=(\eta_1,\ldots,\eta_s) $ be a vector of partitions of the
integer $ d $.  We define the \textit{Hurwitz number}:
$$
    \begin{array}{ccc}
        H^{r}_{h\to g, d}(\underline{\eta}) & := &\mbox{weighted number of}\ \left\{%
        \begin{array}{c}
            degree \ d \ covers \\
            X \stackrel{f}{\longrightarrow} Y\ such\ that:  \\
            \bullet\ $X$\ is \ connected \ of\ genus\ $h$;\ \ \ \ \ \ \ \ \ \ \ \\
            \bullet\ f \ is \ unramified\ over\ \ \ \ \ \ \ \ \ \ \ \ \ \ \ \ \ \ \\
             X\setminus\{p_1,\ldots, p_r,q_1,\ldots,q_s\};\\
            \bullet \ f \ ramifies\ with\ profile\ \eta_i\ over\ q_i; \\ 
            \bullet\ f\ has\ simple\ ramification\ over\ p_i;\\
            \circ\ preimages\ of\ each\ q_i\ with\ same\ \  \ \  \\ ramification\ are\ distinguished\ by\\ appropriate\ markings.  
        \end{array}
        \right\}
    \end{array}
$$
Each cover is weighted by the number of its automorphisms. Figure \ref{Hurw} illustrates the features of this definition.
\end{definition}

\noindent {\bf Remarks.}
\begin{enumerate}
\item For a Hurwitz number to be nonzero, $r$,$g$, $h$ and $\underline{\eta}$ must satisfy the Riemann Hurwitz formula. The above notation is always redundant, and it is common practice to omit appropriate unnecessary invariants.
\item The last condition $\circ$ was  introduced in \cite{gjv:ttgodhn}, and it is well tuned to the applications we have in mind. These Hurwitz numbers differ by a factor of $\prod |\aut(\eta_i)|$ from the classically defined ones where such condition is omitted.
\item One might want to drop the condition of $X$ being connected, and count covers with disconnected domain. Such Hurwitz numbers are denoted by $H^\bullet$.
\end{enumerate}

\begin{example}\label{ex-hn}
\end{example}
\begin{itemize}
	\item $$H^0_{0\to 0, d}((d),(d))= \frac{1}{d}$$
	\item $$H^4_{1\to 0, 2}= \frac{1}{2}$$
	\item $$H^3_{1\to 0, 2}((2),(1,1))= 1$$
\end{itemize}


\subsection{Representation Theory}
 
The problem of computing Hurwitz numbers is in fact a discrete problem and it can be approached using the representation theory of the symmetric group. A standard reference here is \cite{fh:rt}. 
 
\begin{figure}[b]
	\centering
		\includegraphics[width=0.90\textwidth]{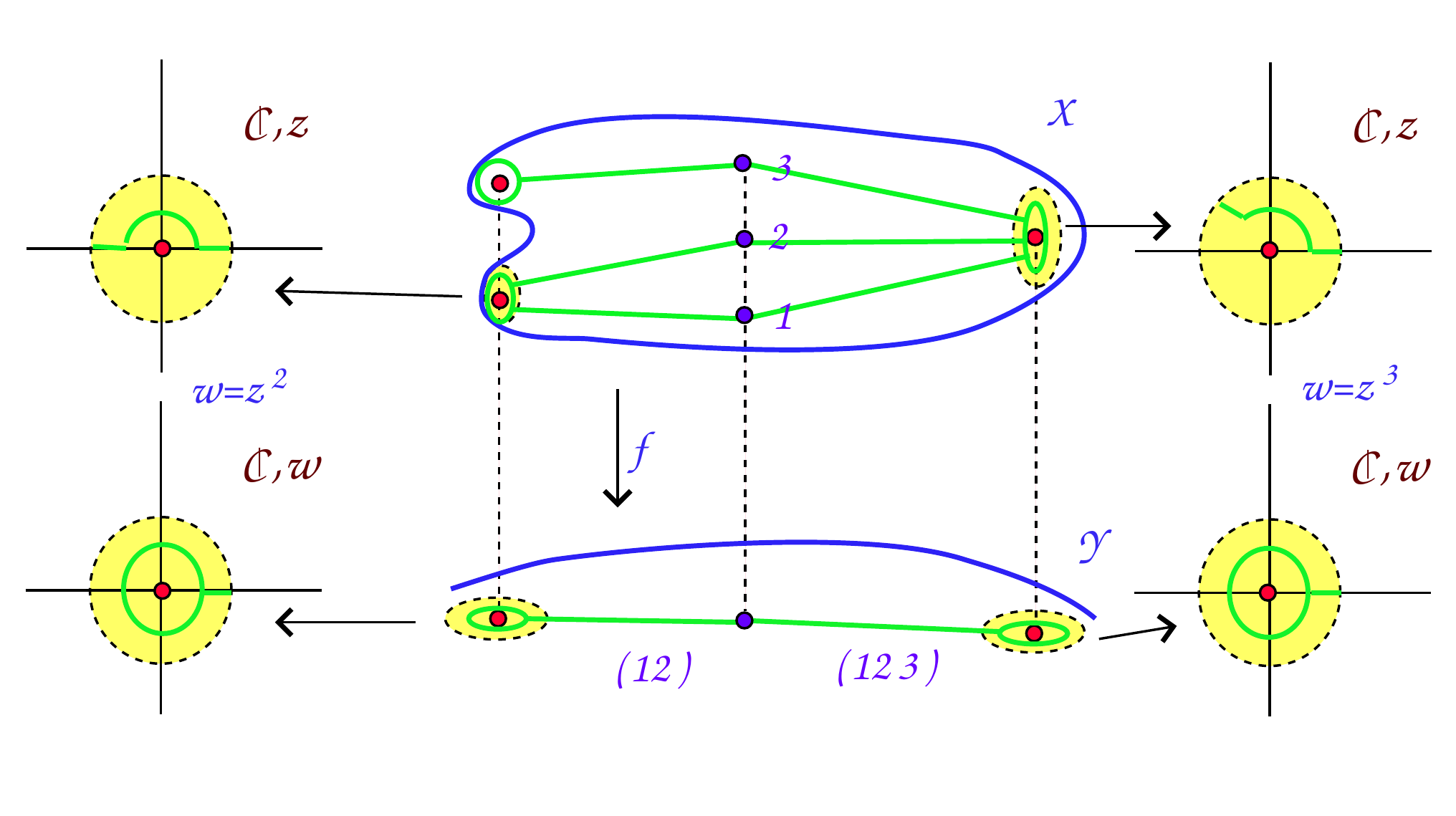}
	\caption{Sketch of the construction of the monodromy representation for the cover $f$ of degree $3$. The base point, at the center, has three labeled inverse images. On the left there is a simple branch point, on the right a branch point with a unique inverse image or ramification index three. Small loops around the branch points and their liftings are depicted. In the ``zoomed areas" it is shown how the lifting of a full turn downstairs results in only part of a turn upstairs: consequently, the liftings of the loops on the base curves are paths and not necessarily loops. }
	\label{fig:monodromy}
\end{figure}

Given a branched cover $f:X \to Y$, pick a point $y_0$ not in the branch locus, and label the preimages $1, \ldots, d$. Then one can naturally define a group homomorphism:
$$
\begin{array}{cccc}
\varphi_f: & \pi_1(Y \minus B,y_0) & \to & S_d\\
           & \gamma            & \mapsto & \sigma_\gamma:\{i\mapsto \tilde{\gamma_i}(1)\},
\end{array}
$$  

where $\tilde{\gamma_i}$ is the lift of $\gamma$ starting at $i$ ($\tilde{\gamma_i}(0)=i$). This homomorphism is called the {\bf monodromy representation}, and its construction is illustrated in Figure \ref{fig:monodromy}.

\noindent{\bf Remarks.}
\begin{enumerate}
	\item A different choice of labelling of the preimages of $y_0$ corresponds to composing $\varphi_f$ with an inner automorphism of $S_d$.
	\item If $\rho\in \pi_1(Y \minus B,y_0)$ is a little loop winding once around a branch point with profile $\eta$, then $\sigma_\rho$ is a permutation of cycle type $\eta$. 
\end{enumerate}

Viceversa, the monodromy representation contains enough information to recover the topological cover of $Y\minus B$, and therefore, by the Riemann existence theorem, the map of Riemann surfaces. To count covers we can count instead (equivalence classes of) monodromy representations. This leads us to the second definition of Hurwitz numbers.

\begin{definition}[Representation Theory]\label{hnrept}
Let $ (Y,p_1,\ldots,p_r,q_1,\ldots,q_s) $ be an $(r+s)$-marked smooth Riemann Surface of genus $g$, and $
\underline{\eta}=(\eta_1,\ldots,\eta_s) $  a vector of partitions of the
integer $ d $:
\begin{equation}
\label{hnrep}
H^{r}_{h\to g, d}(\underline{\eta}):= \frac{|\{\mbox{$\underline{\eta}$-monodromy representations}\ \varphi^{\underline{\eta}}\}|}{|S_d|}\prod |\aut(\eta_i)|,
\end{equation}
where an $\underline{\eta}$-monodromy representation is a group homomorphism
$$
\varphi^{\underline{\eta}}: \pi_1(Y \minus B,y_0)  \to  S_d
$$
such that:\\
 $\bullet$ for $\rho_{q_i}$ a little loop winding around $q_i$ once, $\varphi^{\underline{\eta}}(\rho_{q_i})$ has cycle type $\eta_i$.\\
 $\bullet$ for $\rho_{p_i}$ a little loop winding around $p_i$ once, $\varphi^{\underline{\eta}}(\rho_{p_i})$ is a transposition.\\
 $\star$ Im$(\varphi^{\underline{\eta}}(\rho_{q_i}))$ acts transitively on the set $\{1,\ldots, d\}$.
 \end{definition}
 
\noindent{\bf Remarks.}
\begin{enumerate}
	\item To count disconnected Hurwitz numbers just remove the last condition $\star$.
	\item Dividing by $d!$ accounts simoultaneously for automorphisms of the covers and the possible relabellings of the preimages of $y_0$.
\item $\prod |\aut(\eta_i)|$ is non-classical and it corresponds to condition $\circ$ in Definition \ref{hngeom}.
\end{enumerate}

The count of Hurwitz numbers as in Definition \ref{hnrept} is then translated into a multiplication problem in the class algebra of the symmetric group. 
Recall that $\mathcal{Z}(\C[S_d])$ is a vector space of dimension equal the number of partitions of $d$, with a natural basis indexed by conjugacy classes of permutations.
$$
\mathcal{Z}(\C[S_d])=\bigoplus_{\eta\vdash d} \C \ C_\eta,
$$ 

where 
$$
C_\eta= \sum_{\sigma \in S_d \ \mbox{\tiny of cycle type}\ \eta }{\sigma}. 
$$
We use $|C_\eta|$ to denote the number of permutations of cycle type $\eta$, and $|\xi(\eta)|$ for the order of the centralizer of any permutation of cycle type $\eta$. We also use the notation $C_{Id}= Id$ and $C_\tau= C_{(2,1^{d-2})}$. Define by $\mathfrak{K}$ the element:
\begin{equation}
\mathfrak{K}:= \sum_{\eta \vdash d} |\xi(\eta)| C_\eta^2 \in \cZ\bC[S_d].
\end{equation}

Then the disconnected Hurwitz number is the coefficient of the identity in a product of elements of the class algebra:
\begin{equation} \label{eq:classhurwg}
H^{r\bullet}_{h\to g,d}(\eta_1,\ldots,\eta_s) = \frac{\prod |\aut(\eta_i)|}{d!}[C_e]\mathfrak{K}^g C_\tau^rC_{\eta_s}\ldots C_{\eta_2}  C_{\eta_1} ,
\end{equation}

It is a classical fact that $\mathcal{Z}(\C[S_d])$ is a semisimple algebra with semisimple basis indexed by irreducible representations of $S_d$, and change of bases essentially given by the character table:
\begin{equation}\label{reptoconj}
\mbox{e}_\rho = \frac{\mbox{dim}\rho}{d!}\sum_{\eta\vdash d} \mathcal{X}_{\rho}(\eta) C_\eta
\end{equation}
and 
\begin{equation}\label{conjtorep}
C_\eta = |C_\eta|\sum_{\rho\ \mbox{\tiny irrep. of}\ S_d} \frac{\mathcal{X}_\rho(\eta)}{\mbox{dim}\rho} \mbox{e}_\rho
\end{equation}

By changing basis to the semisimple basis, executing the product there (exploiting the idempotency of the basis vectors) and then changing basis back to isolate the coefficient of the identity, one obtains a closed formula for Hurwitz number - in terms of characters of the symmetric group. This formula is often refereed to as Burnside's Character Formula:
 
\begin{equation}\label{eq:bf}
H^{r,\bullet}_{h\to g,d}(\eta_1,\ldots,\eta_s) = \prod \aut{\eta_i} \sum_{\rho} \left(\frac{\dim  \rho}{d!}\right)^{2-2g} 
\left(\frac{ |C_{\tau}|\chi_{\rho}(\tau)}{\dim  \rho} \right)^r
\prod_{j=1}^s \frac{ |C_{\eta_j}|\chi_{\rho}(\eta_j)}{\dim  \rho}
\end{equation}

\begin{fact}
All Hurwitz numbers are obtained recursively/combinatorially from Hurwitz numbers where the base curve has genus $0$ and there are two  or three branch points. This is a consequence of the {\it degeneration formulas}, that relate the number of covers of a smooth curve with the number of covers of a nodal degeneration of it, and next to the number of cover of its normalization (See  \cite[\S 7.5]{cmhurwitz}).  
For the readers familiar with this language, the above statement is summarized by saying that  Hurwitz numbers form a TQFT (\cite{j:tqft}), whose associated Frobenius Algebra is $\cZ\bC[S_d]$ .
\end{fact}

The name {\bf single (simple) Hurwitz number} (denoted $H^r_g(\eta)$) is reserved for connected Hurwitz numbers to a base curve of genus $0$ and with only one special point where arbitrary ramification is assigned. Recall we adopt a convention introduced by Goulden and Jackson, to mark all ramification points lying over the branch point with special profile. 

The number of simple ramification/branching points, determined by the Riemann-Hurwitz formula, is
\begin{equation}
r= 2g+d-2 +\ell(\eta).
\end{equation}

Single Hurwitz numbers count the number of ways to factor a (fixed) permutation $\sigma\in C_\eta$ into $r$ transpositions that generate $S_d$:
\begin{equation}\label{reduced} 
H^r_g(\eta)= \frac{1}{\prod{\eta_i}}|\{(\tau_1,\ldots,\tau_r\ s.t.\ \tau_1\cdot\ldots\cdot\tau_r=\sigma\in C_\eta, \langle\tau_1,\ldots,\tau_r\rangle=S_d)\}| 
\end{equation}


The first formula for single Hurwitz number was given and ``sort of" proven by Hurwitz in 1891 (\cite{hurwitz}):
$$
H^{r}_0(\eta)= r! d^{\ell(\eta)-3}\prod\frac{\eta_i^{\eta^i}}{\eta_i!}.
$$

Particular cases of this formula were proven throughout the last century, and finally the formula became a theorem in 1997 (\cite{gj:tf}). In studying the problem for higher genus, Goulden and Jackson made the following conjecture.

\begin{conjecture}
For any fixed values of $g,n:=\ell(\eta)$: 
\begin{equation}
\label{gjconj}
H_g(\eta)= r! \prod\frac{\eta_i^{\eta^i}}{\eta_i!} P_{g,n}(\eta_1\ldots,\eta_n),
\end{equation}
where $P_{g,n}$ is a symmetric polynomial in the $\eta_i$'s with:
\begin{itemize}
	\item deg $P_{g,n}= 3g-3+n$;
	\item $P_{g,n}$ doesn't have any term of degree less than $2g-3+n$;
	\item the sign of the coefficient of a monomial of degree $d$ is $(-1)^{d-(3g+n-3)}$. 
\end{itemize}
\end{conjecture}

In \cite{elsv:hnaiomsoc} Ekedahl, Lando, Shapiro and Vainshtein prove this formula by establishing a remarkable connection between simple Hurwitz numbers and tautological intersections on the moduli space of curves (\cite{hm:moc}).

\begin{theorem}[ELSV formula]
For all values of $g,n=\ell(\eta)$ for which the moduli space $\overline{\M}_{g,n}$ exists: 
\begin{equation}
\label{elsv}
H_g(\eta)= r! \prod\frac{\eta_i^{\eta^i}}{\eta_i!} \int_{\overline{\M}_{g,n}}\frac{1-\lambda_1+\ldots+(-1)^g\lambda_g}{\prod(1-\eta_i\psi_i)},
\end{equation}
\end{theorem}

\begin{remark}
In formula \eqref{elsv}, $\psi_i$ is the first Chern class of the $i$-th tautological cotangent line bundle, \cite{k:pc}, and $\lambda_j$ is the $j$-th Chern class of the Hodge bundle, \cite{m:taegotmsoc}. Although these are quite sophisticated geometric objects, in order to understand how the ELSV formula proves Goulden and Jackson's polynomiality conjecture, one only needs to know that $\overline{\M}_{g,n}$ is a compact, smooth orbifold of dimension $3g-3+n$,  $\psi$ classes have degree one, and $\lambda_j$ has degree $j$. 
 The coefficients of $P_{g,n}$ are computed as  intersection numbers on $\overline{\M}_{g,n}$. Using  multi-index notation:
$$
P_{g,n}= \sum_{k=0}^g \sum_{|I_k|=3g-3+n-k}(-1)^k \left( \int \lambda_k \psi^{I_k} \right)\eta^{I_k} 
$$
\end{remark}

\begin{remark}
The polynomial $P_{g,n}$ is a generating function for all linear (meaning where each monomial has only one $\lambda$ class) Hodge integrals on $\overline{\M}_{g,n}$, and hence a good understanding of this polynomial can yield results about the intersection theory on the moduli space of curves. The $ELSV$ formula has given rise to several remarkable applications; a good collection can be found in \cite{kaz}. We recall a couple here:
\begin{description}
\item[\cite{op:mm}] Okounkov and Pandharipande use the ELSV formula to give a proof of Witten's conjecture, that an appropriate generating function for the $\psi$ intersections satisfies the KdV hierarchy. The $\psi$ intersections are the coefficients of the leading terms of $P_{g,n}$, and hence can be reached by studying the asymptotics of Hurwitz numbers:
$$
\lim_{N\to \infty} \frac{P_{g,n}(N\eta)}{N^{3g-3+n}}
$$
\item[\cite{gjv:spolgc}] Goulden, Jackson and Vakil  get a handle on the lowest order terms of $P_{g,n}$ to give a new proof of the $\lambda_g$ conjecture:
$$
\int_{\overline{\M}_{g,n}}\lambda_g\psi^I= {{2g-3+n}\choose{I}}\int_{\overline{\M}_{g,1}}\lambda_g\psi_1^{2g-2}
$$
\end{description}
\end{remark}

We sketch a proof of the $ELSV$ formula following \cite{gv:hnavl}. The strategy is to evaluate an appropriate integral via localization.


We denote 
$$
\M:=\overline{\M}_{g}(\proj^1, \eta\infty) 
$$

the moduli space of relative stable maps of degree $d$ to $\proj^1$, with profile $\eta$ over $\infty$. 
The degenerations  included to compactify are twofold:
\begin{itemize}
	\item away from the preimages of $\infty$ we have degenerations of ``stable maps" type: we can have nodes and contracted components for the source curve, and nothing happens to the target $\proj^1$;
	\item when things collide with $\infty$, then the degeneration is of ``admissible cover" type: a new rational component sprouts from $\infty\in \proj^1$, the special point carrying the profile requirement transfers to this component. Over the node we have nodes for the source curve, with maps satisfying the kissing condition.
\end{itemize}
A more formal description of these spaces and a good picture illustrating the degenerations allowed, can be found in \cite{v:mscgw}, Section $4.10$.
The space $\M$ has virtual dimension $r=2g+d+\ell(\eta)-2$ and admits a globally defined branch morphism (\cite{fp:bd})
$$
br: \M \to Sym^r(\proj^1)\cong \proj^r.
$$
\begin{figure}[b]
	\centering
		\includegraphics[width=0.60\textwidth]{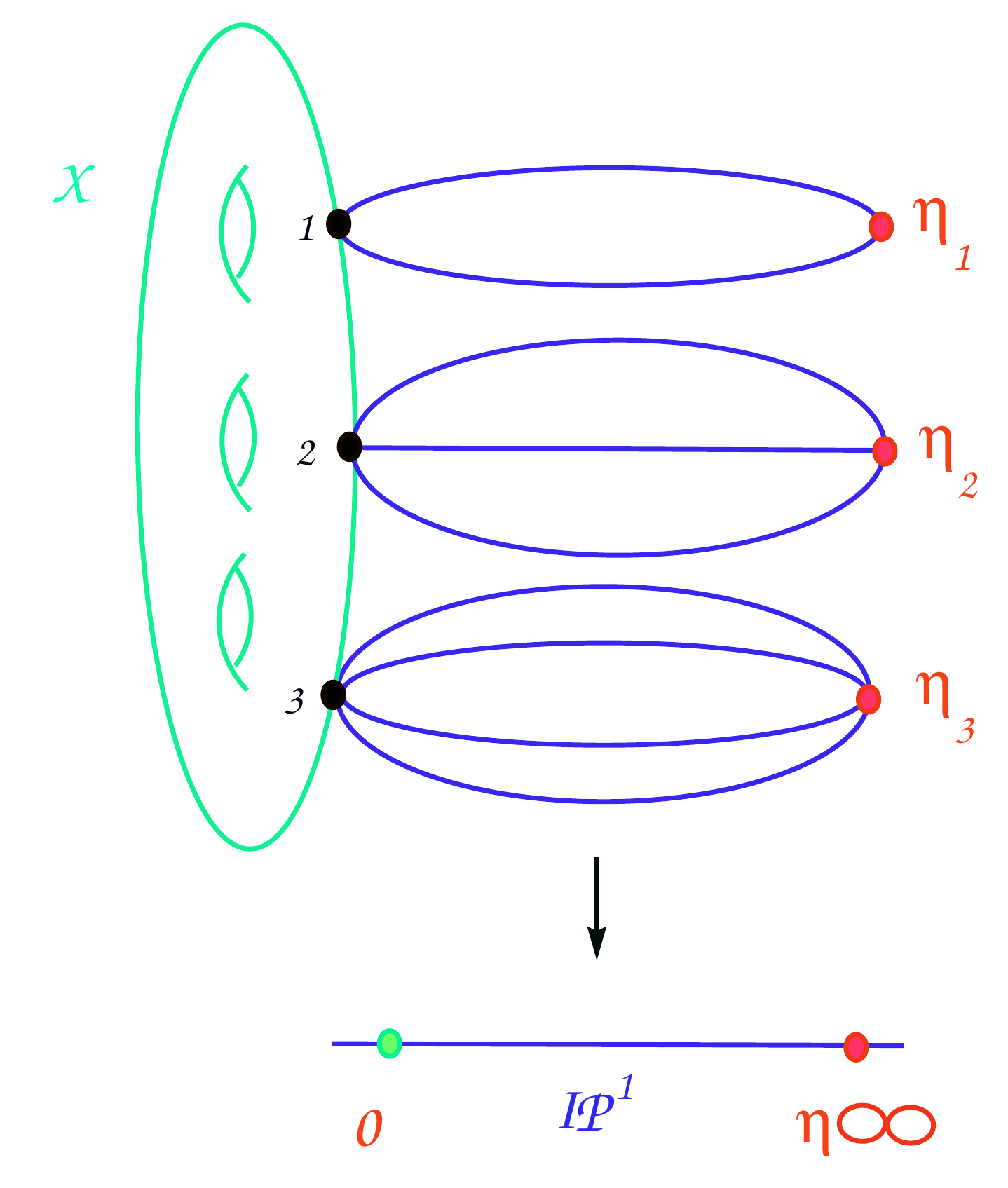}
	\caption{the unique contributing fixed locus in the localization computation proving the $ELSV$ formula. }
	\label{fig:ELSV}
\end{figure}

The simple Hurwitz number:
$$
H^r_g(\eta)=\mbox{deg}(br)= \deg(br^\ast (pt.)\cap [\M]^{vir})
$$
can now interpreted as an intersection number on a moduli space with a torus action and  evaluated via localization. Given a manifold with a torus action, the localization formula allows one to compute the cap product of a cohomology class against the fundamental class in terms of integrals over fixed locus of the torus action, \cite{ab:tmmaec}.  Such technique extends to singular spaces in presence of a virtual fundamental class (\cite{Graber-Pandharipande}). The precise statement of the formula, as well as a hands on introduction to using localization on moduli spaces of stable maps, is found in Chapter 27 of \cite{clay:ms};  virtual localization on spaces of relative stable maps is introduced in \cite{gv:rvl}. The following discussion is intended for readers that are somewhat familiar with this technique.

The map $br$ is made $\C^\ast$ equivariant by inducing the action on $\proj^r$. The key point is to choose the appropriate equivariant lift of the class of a point in $\proj^r$. Recalling that choosing a point in $\proj^r$ is equivalent to fixing a branch divisor, we choose the $\C^\ast$ fixed point corresponding to stacking all ramification over $0$. Then there is a unique fixed locus contributing to the localization formula, depicted in Figure \ref{fig:ELSV}, which is essentially isomorphic to $\overline{\M}_{g,n}$ (up to some automorphism factors coming from the automorphisms of the bubbles over $\proj^1$).

The $ELSV$ formula is obtained from the localization formula. The virtual normal bundle to the unique contributing fixed locus has a denominator part given from the smoothing of the nodes; this produces the denominator with $\psi$ classes in the ELSV formula. One must then compute the restriction to the fixed locus of the equivariant Euler class of $R\pi_\ast f^\ast (T\proj^1(-\infty))$, where $\pi$ is the universal family and $f$ the universal map of the moduli space of relative stable maps. The standard technique for this computation goes as follows: consider the normalization short exact sequence for the generic curve in the fixed locus,  tensor it with $T\proj^1(-\infty)$ and then take the associated long exact sequence in cohomology. Using multiplicativity of the Euler class, one writes the desired class as a product:
\begin{itemize}
\item the Euler class of  a Hodge bundle linearized with weight $1$,  which produces the polynomial in $\lambda$ classes.
\item the Euler class of  a bunch of  trivial but not equivariantly trivial bundles, coming from the covers of the main component. The equivariant Euler class of such bundles is  the product of the corresponding torus weights, and gives rise to the combinatorial pre-factors before the Hodge integral. 
\end{itemize}
For more details on this proof one may consult \cite{gv:taut}, where this localization proof first appeared.

\begin{remark}
An abelian orbifold version of the ELSV formula has been developed by Johnson, Pandharipande and Tseng in \cite{jpt:elsv}. In this case the connection is made  between Hurwitz-Hodge integrals and wreath Hurwitz numbers.
\end{remark}


\section{Double Hurwitz Numbers}
{\bf Double Hurwitz numbers} count covers of $\proj^1$ with special ramification profiles over two points, that for simplicity we assume to be $0$ and $\infty$. Double Hurwitz numbers are classically denoted $H^r_g(\mu,\nu)$; in \cite{cjm:wcfdhn} we start denoting double Hurwitz numbers $H^r_g(\mathbf{x})$, for $\mathbf{x}\in H\subset \R^n$ an integer lattice point on the hyperplane $\sum x_i =0$. The subset of positive coordinates corresponds to the profile over $0$ and the negative coordinates to the profile over $\infty$. We define $\mathbf{x_0}:= \{x_i>0\}$ and $\mathbf{x_\infty}:= \{x_i<0\}$.

The number $r$ of simple ramification is given by the Riemann-Hurwitz formula,
$$
r=2g-2+n
$$
and it is independent of the degree $d$. 
In \cite{gjv:ttgodhn}, Goulden, Jackson and Vakil start a systematic study of double Hurwitz numbers and in particular invite us to consider them as a function:
\begin{equation}
H_g^r (-): \Z^n\cap H \to \Q.
\end{equation}

They prove some remarkable combinatorial property of this function:

\begin{theorem}[\cite{gjv:ttgodhn}] \label{thm-GJVpoly}
The function $H_g(-)$ is a piecewise polynomial function of degree $4g-3+n$.
\end{theorem}

And conjecture some more:
\begin{conjecture}[\cite{gjv:ttgodhn}]
The polynomials describing $H^r_g(-)$ have degree $4g-3+n$, no non-zero coefficients in  degree lower than $2g-3+n$, and are even or odd polynomials (depending on the parity of the leading coefficient).
\end{conjecture}

Shapiro, Shadrin and Vainshtein explore the situation in genus $0$. They describe the chambers of polynomiality by giving the equations of the bounding hyperplanes (:= walls), and give a geometrically suggestive formula for how the polynomials change when going across a wall.

\begin{theorem}[\cite{ssv:g0}]\label{ssvwalls}
The chambers of polynomiality of $H^r_g(-)$ are bounded by \textbf{walls} corresponding to the \textbf{resonance} hyperplanes $W_I$, given by the equation
$$W_I=\left\{\sum_{i\in I} x_i=0 \right\},$$
for any $I\subset\{1,\ldots,n\}$.

Let $\cham_1$ and $\cham_2$ be two chambers adjacent along the wall $W_I$, with $\cham_1$ being the chamber with $x_I<0$.  The Hurwitz number $H^r_g(\mathbf{x})$ is given by polynomials, say $P_1(\mathbf{x})$ and $P_2(\mathbf{x})$, on these two regions.  A wall crossing formula is a formula for the polynomial $$WC_I^r(\mathbf{x})=P_2(\mathbf{x})-P_1(\mathbf{x}).$$

Genus $0$ wall crossing formulas have the following inductive description:
\begin{equation}
\label{g0wc}
WC_I^r(\mathbf{x})= \delta {{r}\choose{r_1,r_2}} H^{r_1}(\mathbf{x_I},\delta) H^{r_2}(\mathbf{x_{I^c}},-\delta),
\end{equation}
where $\delta= \sum_{i\in I}x_i$ is the distance from the wall at the point where we evaluate the wall crossing.
\end{theorem}

\noindent
{\bf Remarks.}
\begin{enumerate}
	\item This formula appears not to depend on the particular choice of chambers $\cham_1$ and $\cham_2$ that border on the wall, but only upon the wall $W_I$; however the polynomials for the simpler Hurwitz numbers in the formula  depend on chambers themselves.
	\item The walls $W_I$ correspond to values of $\mathbf{x}$ where the cover could potentially be disconnected, or where $x_i=0$ for some $i$. In the first case the formula reminds of a boundary divisor degeneration formula, and somehow begs for a geometric understanding.
	\item  Crossing the second type of wall corresponds to moving a ramification between $0$ and $\infty$.  In the traditional view of double Hurwitz numbers, these were viewed as separate problems: the lenghts of the profiles over $0$ and $\infty$ were considered two independent discrete invariants.  However, here we see that the more natural invariant is just the total length of the special ramification imposed: this motivates  $\mathbf{x}$ replacing $\mu,\nu$ in our notation. In genus $0$ the wall crossing formula for ${x}_i=0$ is trivial - and as such identical to all other wall crossing formulas.  In higher genus this second type of wall crossing are not trivial any more, while still obeying the same wall crossing formulas as wall crossing of the first type.
\end{enumerate}

The way Goulden, Jackson and Vakil prove their result is similar to \cite{op:mm}: they compute double Hurwitz numbers by counting decorated ribbon graphs on the source curve. A ribbon graph is  obtained by pulling back a set of arcs from the target sphere (connecting $0$ to the simple ramification points) and then stabylizing. Each ribbon graph comes with combinatorial decorations that are parameterized by integral points in a polytope with linear boundaries in the $x_i$'s. Standard algebraic combinatorial techniques  (Ehrhart theory) then show that such counting yields polynomials so long as the topology of the various polytopes does not change. The downside of this approach is that these are ``large" polytopes  (namely $4g-3+n$ dimensional) and it is hard to control their topology.


The approach of \cite{cjm:thn} to this problem is motivated by tropical geometry. Double Hurwitz numbers are computed in terms of some trivalent polynomially weighted graphs called {\bf monodromy graphs} that are, in a sense, ``movies of the monodromy representation": they encode the cycle type of all successive compositions by transpositions of the initial permutation $\sigma$ giving the monodromy over $0$. Monodromy graphs can be thought as tropical covers, even though this point of view is not necessary other than to give the initial motivation. This combinatorial encoding gives a straightforward and clean proof of the genus $0$ situation, contained in Section 6 of \cite{cjm:thn}. In \cite{cjm:wcfdhn}, we show that in higher genus each graph $\Gamma$ comes together with a polytope $P_\Gamma$ (with homogenous linear boundaries in the $x_i$) and we have to sum the polynomial weight of the graph over the integer lattice points of $P_\Gamma$. It is again standard (think of it as a discretization of integrating a polynomial over a polytope) to show that this contribution is polynomial when the topology of the polytope does not change. The advantage  is that we have transferred most of the combinatorial complexity of the enumerative problem in the polynomial weights of the graph: our polytopes are only $g$ dimensional and it is possible to control their topology. 
The result is a  wall-crossing formula that generalizes \cite{ssv:g0} to arbitrary genus.
\begin{theorem}[ \cite{cjm:wcfdhn}] \label{maintheorem}
\begin{equation}
\label{wcformula}
WC_I^r(\mathbf{x})=\sum_{\substack{s+t+u=r \\|\mathbf{y}|=|\mathbf{z}|=|\mathbf{x}_I|}} (-1)^{t} {r \choose s,t,u}\frac{\prod\mathbf{y}_i}{\ell(\mathbf{y})!}\frac{\prod \mathbf{z}_j}{\ell(\mathbf{z})!} H^{s}(\mathbf{x}_I, \mathbf{y})H^{t\bullet}(-\mathbf{y},\mathbf{z})H^{u}(\mathbf{x}_{I^c},- \mathbf{z})
\end{equation}
Here $\mathbf{y}$ is an ordered tuple of $\ell(\mathbf{y})$ positive integers with sum $|\mathbf{y}|$, and similarly with $\mathbf{z}$.
\end{theorem}

\subsection{ELSV Formula for Double Hurwitz Numbers}
\label{sec:elsvdhn}

The combinatorial structure of double Hurwitz numbers seems to suggest the existence of an $ELSV$ type formula, i.e. an intersection theoretic expression that explains the  piecewise polynomiality properties discussed above.
This proposal was initially made in \cite{gjv:ttgodhn} for the specific case of \textit{one-part} double Hurwitz numbers, where there are no wall-crossing issues. After \cite{cjm:wcfdhn}, Bayer-Cavalieri-Johnson-Markwig propose an intriguing, albeit maybe excessively bold strengthening of Goulden-Jackson-Vakil's original conjecture, which we now report:

\begin{conjecture}
For $\mathbf{x}\in \Z^n$ with $ \sum x_i=0$,
\begin{equation}
H_g(\mathbf{x})= \int_{\overline{P}(\mathbf{x})}\frac{1-\Lambda_2+\ldots+(-1)^g\Lambda_2g}{\prod(1-x_i\psi_i)},
\end{equation}
where
\begin{enumerate}
	\item $\overline{P}(\mathbf{x})$ is a (family of) moduli space(s, depending on $\mathbf{x}$) of dimension $4g-3+n$.
	\item $\overline{P}(\mathbf{x})$ is constant on each chamber of polynomiality.
	\item The parameter space for double Hurwitz numbers can be identified with a space of stability conditions for a moduli functor and the $\overline{P}(\mathbf{x})$ with the corresponding compactifications.
	\item $\Lambda_{2i}$ are tautological Chow classes of degree $2i$.
	\item $\psi_i$'s are cotangent line classes.
\end{enumerate}
\end{conjecture}

\noindent{\bf Remarks on one-part Double Hurwitz Numbers.}
\begin{itemize}
	\item Goulden, Jackson and Vakil restrict their attention to the one part double Hurwitz number case, where there are no issues of piecewise polynomiality. Here they propose that the mystery moduli space may be  some compactification of the universal Picard stack over $\overline{\mathcal{M}}_{g,n}$. They verify that such a conjecture holds for genus $0$\footnote{ The paper \cite{gjv:ttgodhn} claims to verify the conjecture  for genus $1$ by identifying $\overline{Pic}_{1,n}$ with $\overline{\mathcal{M}}_{1,n+1}$, but the referee pointed to my attention that their argument is not, or at least not immediately, correct.}.
	\item Defining the symbols $\langle\langle \vec{\tau_k} \Lambda_{2i} \rangle\rangle$ to be the appropriate coefficients of the one-part double Hurwitz polynomial, they are able to show that such symbols satisfy the {\it string} and {\it dilaton} equations. 
	\item The computations:
	$$\langle\langle \tau_{k_1} \ldots \tau_{k_n} \Lambda_{2g} \rangle\rangle={{2g+n-3}\choose{k_1, \ldots, k_n}} \langle\langle \tau_{2g-2} \Lambda_{2g} \rangle\rangle$$
	$$\langle\langle \tau_{2g-2} \Lambda_{2g} \rangle\rangle= \frac{2^{2g-1}-1}{2^{2g-1}(2g)!}|B_{2g}|= \int_{\overline{M}_{g,1}} \lambda_g \psi^{2g-2}$$
	are used as evidence for the following conjecture
	\begin{conjecture}[\cite{gjv:ttgodhn}]
	 There is a natural structure morphism $\pi: \overline{P}(\mathbf{x}) \to \overline{M}_{g,n}$ and $ \pi_\ast \Lambda_{2g} = \lambda_g$.
	\end{conjecture}
	\item Shadrin-Zvonkine, Shadrin (\cite{sz:cov},\cite{s:ots}) show that the symbols $\langle\langle \tau_{k_1} \ldots \tau_{k_n} \Lambda_{2i} \rangle\rangle$ can be organized in an appropriate generating function to satisfy the Hirota equations (this implies being a $\tau$ function for the KP hierarchy - in the pure ``descendant" case this is an analog of Witten's conjecture for the moduli space of curves).
\end{itemize}
\noindent{\bf Remarks on Wall-Crossings.}
\begin{itemize}
\item An observation of A. Craw: the polytopes appearing in the proof of the wall-crossing formula of \cite{cjm:wcfdhn} have constant topology precisely on chambers of $\theta$-stabiity of the quiver varieties associated to the graph giving rise to the polytope.
\item Kass and Pagani have been studying a variation of the $\Theta$ divisor on a family of compactified Jacobians $J_{g,n}(\phi)$ (where $\phi$ varies in an affine space $V_{g,n}$  and, in vague terms, keeps track of the multidegree of a line bundle on a reducible curve), which involves boundary corrections as walls are crossed. It is very plausible that such classes might appear in some kind of intersection theoretic formula for the double Hurwitz numbers on the universal Picard, and contribute to the wall-crossing phenomenon.
 \end{itemize}

\section{Geometry behind Double Hurwitz Numbers} \label{steffen}
This section contains two alternative geometric approaches to the study of double Hurwitz numbers. First we investigate the wall crossing phenomena as a ``variation" of $\psi$ classes on the moduli space of curves. Then we present the point of view offered by tropical geometry.

\subsection{An Intersection Theoretic Formula on $\overline{M}_{g,n}$}
\begin{figure}[tb]
\includegraphics[width=0.8\textwidth]{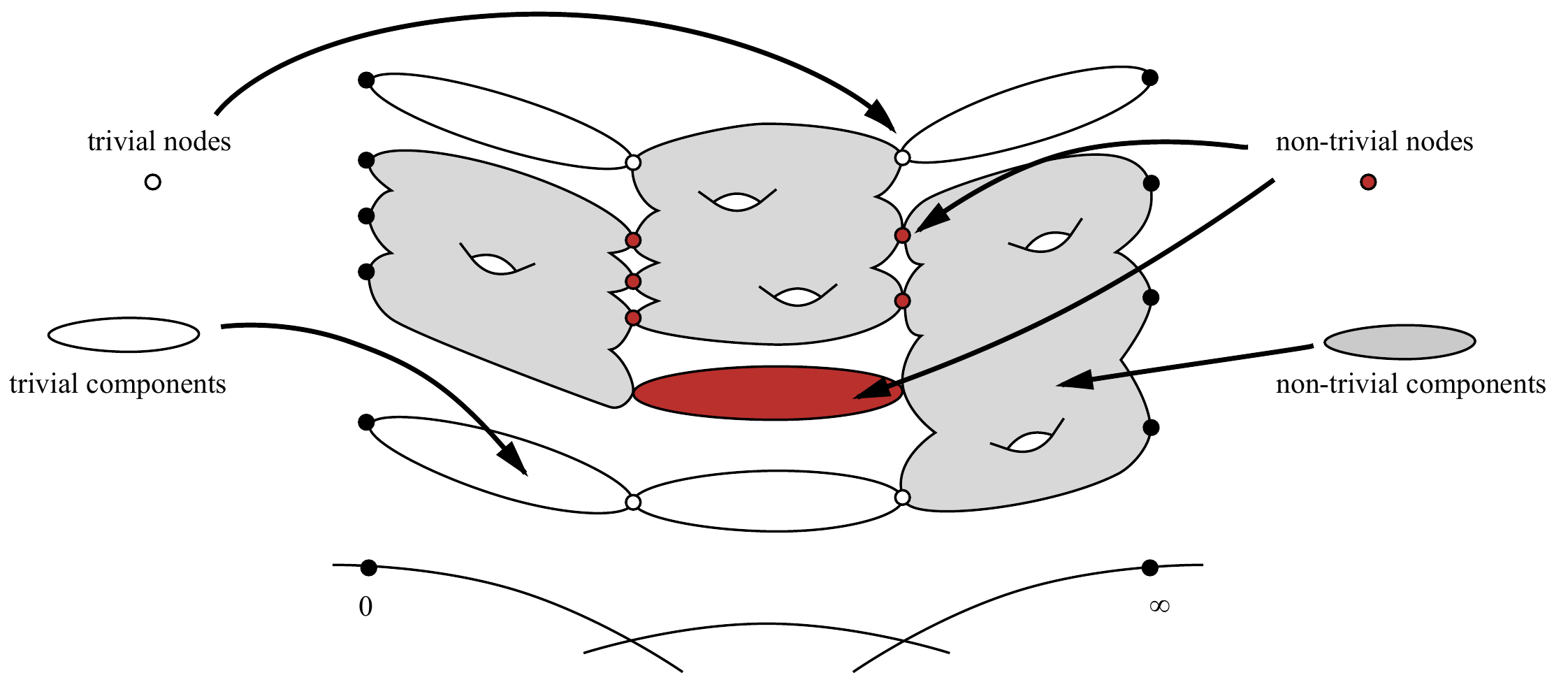}
\caption{ A map in the boundary of the space of rubber relative maps.  Trivial components of the source curve are rational components which becomes unstable when the map is forgotten. Nodes that attach to at least one trivial component are called trivial nodes. A convenient abuse of notation: a trivial component that, when contracted, gives rise to a non-trivial node, is called  a non-trivial node.}
\label{rubmap}
\end{figure}

An alternative approach giving  a more geometric view of the structure of double Hurwitz numbers is pursued in \cite{cm:dhn}. 
We must work with a variant to the moduli space of relative stable maps, where two maps are identified if they differ by an automorphism of $\PP^1$ that preserves $0$ and $\infty$, i.e. $f \equiv \alpha f$ for any $\alpha\in \bC^\ast$. Since "stretching" $\PP^1$ by elements of its torus respects these equivalence classes, the corresponding moduli space was named space of maps to a rubber $\PP^1$ and then, by abuse of notation, it became just ``rubber maps".
Denote by $\relrub{g}{\bx} $ the moduli space of rubber relative stable maps to $\PP^1$, relative to two points with ramification profiles determined by $\mathbf{x}$.  There is a natural forgetful morphism $stab: \relrub{g}{\bx} \to \cms_{g,n} $, which  remembers the (stabilization of the) source curve of the maps. The advantage of working with rubber maps is that $stab$ does not have positive dimensional fibers. Figure \ref{rubmap} illustrates a possible rubber map and introduces some relevant notation.

\begin{figure}[tb]
\includegraphics[width=0.9\textwidth]{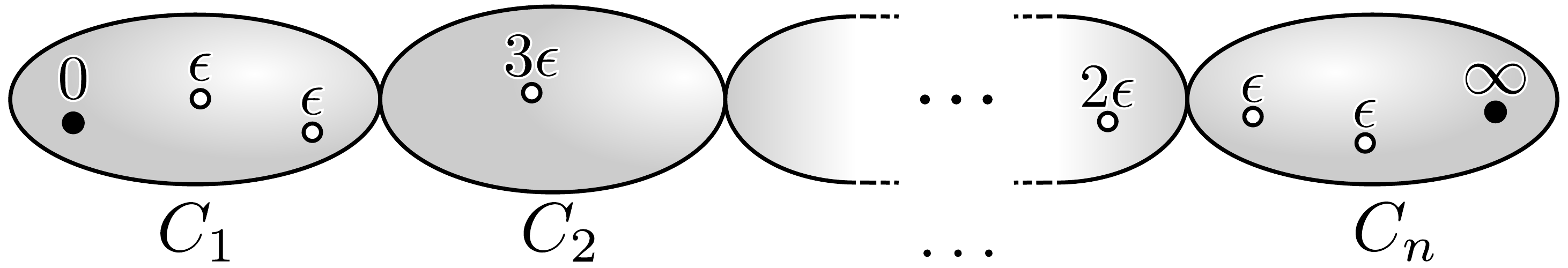}
\caption{The moduli space $\cms_{br}$ parameterizes configurations of $r+2$ marked points on chains of rational curves. The two special points are labeled $0$ and $\infty$ and must be smooth points one on each terminal component of the chain. The remaining $r$ points are undistinguishable, hence in the picture we label them all $\epsilon$; they may not coincide with nodes or with $0$ and $\infty$, but they may coincide among themselves: in the picture we indicate this with multiplicities - a point labeled $3\epsilon$ means that three points have come together. A configuration is stable if every component of the chain has at least three distinct special points (marks or nodes). Two configurations are equivalent if they differ by an automorphism of the chain. }
\label{branch}
\end{figure}

There is also a natural branch morphism to a moduli space we denote $\cms_{br}$, which is described in Figure \ref{branch}.
 For readers familiar with  this language, $\cms_{br}$ is the stack quotient of a Losev-Manin space (which we interpret as an appropriate Hasset space of weighted pointed rational curves) by the symmetric group permuting the ``light points". 
 The diagram of spaces and natural morphisms
\begin{equation}\label{diag:central}
\xymatrix{  \relrub{g}{\bx} \ar[r]^{stab}  \ar[d]_{br} & \cms_{g,n} \\
\cms_{br} :=\left[\cms_{0,2+r}(1,1,\epsilon,\ldots,\epsilon)/S_r\right]
}
\end{equation}
allows to express the double Hurwitz number as the degree of a  cycle on $\cms_{g,n}$:
\beq
\label{eq:dhnif}
H^r_g(\bx)[pt.]= \stab_\ast(\br^\ast([pt.]) \cap [ \relrub{g}{\bx} ]^{vir}).
\eeq

We rewrite \eqref{eq:dhnif} in terms of $\psi$ classes. We have three different kinds of $\psi$ classes:

\begin{enumerate}
	\item $\hat{\psi}_0$:  the psi class on the universal target space at the relative divisor $0$, i.e. the first Chern class of the cotangent line bundle at the relative point $0$. 
		\item $\tilde{\psi}_i$: the psi classes on the space of rubber stable maps at the $i$-th mark on the source curve. Remember that we are marking the preimages of the relative divisors.
	\item $\psi_i$ is the ordinary psi class on the moduli space of curves.	 
\end{enumerate}

These different $\psi$ classes are related via the tautological morphisms as follows:


\beq \label{eq:psi2}
br^\ast(\hat\psi_0)= x_i\tilde{\psi_i}
\eeq
whenever $x_i>0$ (i.e. the $i$-th mark is a preimage of $0$). Denote by $\sD_i^\fc$  the divisor parameterizing maps where the $i$-th mark is supported on a trivial component. Then:

\begin{equation} \label{eq:psi3}
\widetilde{\psi_i} = stab^\ast\psi_i+ \frac{1}{x_i}\sD^\fc_i.\end{equation}


We note that we can express the class of a point on $\cms_{br}$  as a multiple of the top power of $\hat{\psi}_0$:
\beq \label{eq:psi1}
\hat\psi_0^{r-1}= \frac{1}{r!}[pt.]
\eeq
Combining equations \eqref{eq:psi2},\eqref{eq:psi3} and \eqref{eq:psi1} provides the chain of equalities:
\begin{align*}
\br^\ast[\pt]\cap\left[\relrub{g}{\bx}\right]^{\rm vir} \stackrel{\eqref{eq:psi1}}{=}& r!\br^\ast\left(\widehat{\psi}_0^{r-1}\right)\cap\left[\relrub{g}{\bx}\right]^{\rm vir}\\
\stackrel{\eqref{eq:psi2}}{=}& r!\left(x_i\widetilde{\psi}_i\right)^{r-1}\cap\left[\relrub{g}{\bx}\right]^{\rm vir}\\
\stackrel{\eqref{eq:psi3}}{=}& r!\left(x_i\stab^\ast\psi_i+\sD_i^{\fc}\right)^{r-1}\cap\left[\relrub{g}{\bx}\right]^{\rm vir}
\end{align*}
in the Chow ring of $\relrub{g}{\bx}$ for any $i$-th marked pre-image of $0$.  Pushing forward via the stabilization morphism, expanding, and applying the projection formula yields:
\begin{align}
H_g(\bx)[\pt] =& r! \stab_\ast\left(\left(x_i\stab^\ast\psi_i+\sD_i^{\fc}\right)^{r-1}\cap\left[\relrub{g}{\bx}\right]^{\rm vir}\right) \nonumber\\
=& r! \stab_\ast\left( \sum_{k=0}^{r-1}\binom{r-1}{k}(x_i\stab^\ast\psi_i)^{r-1-k}\left(\sD_i^{\fc}\right)^{k}\cap\left[\relrub{g}{\bx}\right]^{\rm vir}\right)\nonumber\\
=& r! \sum_{k=0}^{r-1}\binom{r-1}{k}(x_i\psi_i)^{r-1-k}\stab_\ast\left( \left(\sD_i^{\fc}\right)^{k}\cap\left[\relrub{g}{\bx}\right]^{\rm vir}\right). \label{nice}
\end{align}

It is now a matter of bookkeeping to show that formula \eqref{nice} proves that $H_g(\bx)$ is a piecewise polynomial function of degree $4g-3+n$. The  required conceptual ingredients  are:
\begin{itemize}
\item The supports of the classes $\stab_\ast\left(( \sD_i^{\fc})^k\right)$ are constant in the chambers of polynomiality of the double Hurwitz number.
\item The push-forward of a stratum in the moduli space of relative stable maps contains a ``ghost automorphism" factor for each non-trivial node in the generic (source) curve parameterized by the stratum. Thus each non-trivial node contributes a linear factor in $x_i's$ and in $g$ additional variables $l_i$.
\item The polynomial function given by the product of the node factors must be summed over the lattice points of a $g$-dimensional polytope whose boundaries are linear in the $x_i$'s  and $l_i$'s.
\item The push-forward of the virtual class of the moduli space of relative stable maps intersects monomials in $\psi$ classes as if it were an even polynomial class of degree $2g$ (\cite{bssz:psi}).
\end{itemize} 
The class $\stab_\ast(\left[\relrub{g}{\bx}\right]^{\rm vir})$ is called the {\bf double ramification cycle} and it will be the hero of Section \ref{sec:drc}.

\subsection{Tropical Point of View}
\label{trop}
In overly simplistic terms, tropical geometry replaces classical geometric objects with corresponding piecewise linear objects. Tropical invariants are then combinatorial quantities that are accessible for computations. In many circumstances, despite the severe degeneration of the objects and apparent loss of information in the process, tropical invariants end up recovering the corresponding classical ones. While in the early days of tropical geometry these correspondences had a bit of a {\it magical} feel to it, we are now coming to a more conceptual understanding of why tropical geometry ``works so well".

In \cite{cjm:thn}, we defined {\bf tropical Hurwitz numbers} as follows:
\begin{enumerate}
\item  constructed a polyhedral complex parameterizing  affine linear maps of decorated metric graphs (tropical covers).
\item defined a natural {\it branch morphism} recording the images of the vertices of the source graphs, together with natural linear algebraic multiplicities associated to the restriction of the morphism to maximal cells in the polyhedral complex.
\item defined tropical Hurwitz numbers to be the degree of the branch morphism.  
\end{enumerate}

The correspondence theorem stemmed from the fact that the computation of the degree of the tropical branch morphism is  a sum over graphs which  parallels a classical Hurwitz count through monodromy representations. Each graph corresponds to a sequence of cycle types of all successive products of generators of the fundamental group of the base. The tropical weights just {\it happen} to be the multiplicities coming from the cut and join equation for multiplication by a transposition in the class algebra of $S_d$. In a similar fashion, 
Bertrand-Brugalle-Mikhalkin \cite{BBM} prove a correspondence theorem for arbitrary Hurwitz numbers by observing that the tropical counts organize the classical computation of Hurwitz numbers via a maximal degeneration of the target.

In all cases, a direct link between the geometry of the classical and tropical objects to be counted (and the corresponding moduli spaces) was lacking!

The major step in filling this gap is made in work of Abramovich-Caporaso-Payne \cite{ACP}, who create a connection between the classical and tropical moduli spaces of curves by ``going through" non-archimedean geometry. More specifically, they consider  the Berkovich analytification $\overline{\sM}_{g,n}^{an}$, consisting of points over a valued field together with an extension of the field valuation.  The resulting space is quite complicated to describe, but it is a compact metric space with a Hausdorff topology. The really nice feature however arises once we consider the a toroidal structure of $\overline{\sM}_{g,n}$ given by its modular boundary. Roughly speaking, a divisor $D$ on a space $X$ induces a toroidal structure if, locally in some appropriate topology, $D$ looks like the boundary of some toric variety. Since the boundary of $\overline{\sM}_{g,n}$ is normal crossings, this condition is certainly met. This natural toroidal structure of the moduli space of curves identifies an (extended) cone complex onto which the Berkovich space retracts, called {\bf skeleton} $\Sigma\overline{\sM}_{g,n}^{an}$. Cones in this complex are in bijection with irreducible boundary strata in the moduli space of curves, that are identified by the dual graph of the general curve in the stratum. Let $C_\Gamma$ be the cone corresponding to the dual graph $\Gamma$. Roughly speaking, one should think of a point if $C_\Gamma$ as  an infinitesimal family of curves such that the central fiber has dual graph equal to $\Gamma$ and the generic fiber may be smooth (but it doesn't have to be).

There is a natural tropicalization morphism  between $\Sigma\overline{\sM}_{g,n}^{an}$  and $\overline{\sM}_{g,n}^{trop}$, the parameter space of metric graphs of genus $g$ with $n$ ends. To a point in the skeleton one associates the dual graph to the central fiber. In order  to give a metrization of an edge of the graph, one looks at the local analytic equation $xy=t$ of the family in a neighborhood of the corresponding node. Recalling that $t$ is an element of a field with a valuation, it is natural to assign to the corresponding edge length equal to the valuation of $t$. 

\begin{theorem}[\cite{ACP}] 
The tropicalization morphism
$$
trop: \Sigma\overline{\sM}_{g,n}^{an} \to \overline{\sM}_{g,n}^{trop}
$$
is an isomorphism of (generalized, extended) cone complexes.
\end{theorem}

In recent work with Markwig and Ranganathan \cite{CMR14}, we apply the techniques of \cite{ACP} to the setting of admissible covers. The (normalization of the) moduli space of admissible covers, which is an (orbifold smooth) compactification of the space of maps from curves to $\PP^1$ with specified ramification data, also has normal crossings boundary, and hence a natural toroidal structure. Hence we have a skeleton inside the analytification  of the space of admissible covers. There is also a tropicalization map that assigns to an infinitesimal family of admissible covers the pair of dual graphs of the source and target curves in the central fiber with edges metrized by the valuation of the corresponding node smoothing parameters. Such decorated graphs are defined to be {\it tropical admissible covers} in \cite{gonality}.
We organize tropical admissible covers into a cone complex that we call the  moduli space of tropical admissible covers and prove a result similar in flavor to \cite{ACP}. 
 
\begin{theorem}[\cite{CMR14}]~\label{thm: theorem1}
The tropicalization map $$trop: \Sigma \cH^{an}_{g\to h, d}(\vec\mu)\to \cH^{trop}_{g\to h, d}(\vec\mu)$$  is a face morphism of cone complexes, i.e. the restriction of $trop$ to any cone of $\Sigma(\Hbar^{an}_{g\to h,d}(\vec\mu))$ is an isomorphism onto a cone of the tropical moduli space $\cH^{trop}_{g\to h,d}(\vec\mu)$. 

Let $br$ denote the branch map $\Hbar^{}_{g\to h, d}(\vec \mu)\to \Mbar_{h,r+s}$. Then the following diagram is commutative:

\begin{equation}
\xymatrix{
\Hbar^{an}_{g\to h, d}(\vec \mu) \arrow[r]_{trop} \arrow[d]_{br^{an}}   & \Hbar^{trop}_{g,d}(\vec \mu) \arrow[d]_{br^{trop}} \\
\Mbar^{an}_{h,r+s} \arrow[r]^{trop} & \Mbar^{trop}_{h,r+s}, & 
}
\label{diagbr}
\end{equation}
\end{theorem}

 In the sense of Abramovich-Caporaso-Payne, the tropical moduli space of admissible covers is a tropicalization of the classical one. Our result states that the tropicalization morphism maps cones isomorphically to cones but is typically  not globally one to one. The morphism is surjective\footnote{Here there is a slight issue of conventions in the literature. In \cite{gonality}, a tropical admissible cover may not represent the dual graph of a classical admissible cover, due to a local issue called the {\it Hurwitz existence problem}. In our work, we find it natural to include the positive solution to the Hurwitz existence problem into our definition of tropical admissible cover.} but not injective, due to the fact that distinct analytic covers may have the same dual graphs, and to multiplicities arising when normalizing the space of admissible covers. A little bit of work goes into keeping track of the degree of the fibers, but once that is done, we can extract enumerative information from the commutativity of $\eqref{diagbr}$.
 We  compare the degree of the analytic branch map and the degree of the tropical branch map and  recover a purely geometric proof of the correspondence theorem of Hurwitz numbers of \cite{cjm:thn} and \cite{BBM}.

\section{The Double Ramification Cycle}
\label{sec:drc}

Let $(C, p_1, \ldots, p_n)$ be a smooth genus $g$,  $n$-pointed curve, and let $\mathbf{x}\in \bZ^n$ be a vector of non-zero integers adding to $0$. We state  three equivalent conditions:
\begin{description}
	\item[PD] The divisor $D(\mathbf{x})=\sum x_ip_i$ is principal.
	\item[M] There exists a map $f: C\to \PP^1$ such that the ramification points over two special points of $\PP^1$ are  the $p_i$'s and the ramification profiles are given by the positive and negative parts of $\mathbf{x}$. 
		\item[J]  Consider the section $\Phi_{\mathbf{x}}: \sM_{g,n}\to J_{g,n}$ defined by  $\Phi_{\mathbf{x}}(C,p_1, \ldots, p_n)= (C, \sum x_ip_i)$. Then the image of $\Phi_{\mathbf{x}}(C,p_1, \ldots, p_n)$ lies in the $0$-section of the universal Jacobian.
\end{description}

For a fixed $\mathbf{x}$, we wish to describe a {\bf meaningful} compactification of the locus of curves in $\sM_{g,n}$ satisfying the above conditions.  This goal is often referred to as {\it Eliashberg's Problem}, who posed this question  in the early $2000$'s for the development of symplectic field theory. Such a compactification would define a cycle in the Chow ring of $\overline{\sM}_{g,n}$, which we call the {\bf double ramification cycle} ($DRC(\mathbf{x})$).

The word meaningful is purposefully ambiguous. Here is a few things it could mean:
\begin{enumerate}
	\item $DRC(\mathbf{x})$ has a modular interpretation.
	\item $DRC(\mathbf{x})$ has good geometric properties.
	\item $DRC(\mathbf{x})$ has good algebro-combinatorial properties.
\end{enumerate}

And here are a few candidates for a rigorous definition of the double ramification cycle:
\begin{enumerate}
\item $\stab_\ast([\overline{M}^{\sim}_{g}(\PP^1, \mathbf{x})]^{vir})$: the push-forward of the virtual fundamental class of the moduli space of rubber relative stable maps to $\PP^1$.
\item $\stab_\ast([\overline{Adm}_{g}(\PP^1, \mathbf{x})])$: the push-forward of the fundamental class of the moduli space of admissible covers to $\PP^1$.
\item the pullback of the $0$-section of some compactified universal Jacobian via some resolution of indeterminacies of the section $\sigma$ described above. 
\end{enumerate}

In \cite{H11}, Hain proposes  a class for the double ramification cycle in the tautological ring of the moduli space of curves of compact type:
$$
DRC(\mathbf{x}) = \Phi_{\mathbf{x}}^\ast ([0]),
$$
he points  out that the abel map $\Phi_{\mathbf{x}}$ extends without indeterminacy to curves of compact type. Further, on an abelian variety we have the equation
$$
[0] = \frac{\Theta^{g}}{g!},
$$
 which is then extended to hold in a family over the image of $\Phi_{\mathbf{x}}$.

Then using the theory of normal functions he computes the pull-back of the theta divisor in terms of standard tautological classes in the moduli space of curves, to obtain a class which we call the {\bf Hain class} and denote $\mathfrak{h}_g(\mathbf{x})$:
\beq \label{eq:hain}
\mathfrak{h}_g(\mathbf{x})= \frac{1}{g!}\left( 
-\frac{1}{4} \sum_{h, \mathcal{S}}
\left(
\sum_{i\in \mathcal{S}} x_i
\right)^2 \delta_{h,\mathcal{S}}
\right)^g;
\eeq
 the summation is  for $h=0, \ldots, g$ and over $\mathcal{S}\subseteq [n]$, i.e. all subsets of the set of marks. The symbol $ \delta_{h,\mathcal{S}}$ denotes the class of the boundary divisor parameterizing nodal curves with irreducible components of genera $h, g-h$  and the marks in the subset $\mathcal{S}$ on the component of genus $h$. We also adopt the following conventions for the unstable cases:
$$
\delta_{0, \{i\}}= \delta_{g, [n] -i}= -\psi_i \ \ \ \ \ \ \ \ \ \ \ \ \ \ \ \ \
\delta_{0, \phi}= \delta_{g, [n]}= \kappa_1.
$$ 

\noindent{\bf Remarks.}
\begin{enumerate}
	\item In \cite{GZ1}, Grushevski-Zakharov give a simpler computation of the class $\Phi_{\mathbf{x}}^\ast(\Theta)$ via test curves.
	\item Equation \eqref{eq:hain} exhibits the Hain class as a homogeneous polynomial class of degree $2g$.
	\item In their recent paper \cite{mw:fund}, Marcus and Wise show that the restriction to compact type of the Hain class coincides with the restriction of compact type of the push-forward of the virtual class of the moduli space of rubber relative stable maps.
	\item Although  Hain's formula determines a very nice class in the tautological ring of the full moduli space of curves, it is not a good candidate to be the full $DRC(\mathbf{x})$. For example the appropriate intersection number with a top power of a psi class in the one-part case does not give the double Hurwitz number. We use this simple fact to compute the  correction between the Hain class and the DRC in genus $1$.	
\end{enumerate}
\noindent {\it Genus $1$.}
We restrict our attention to genus $1$ and two marked points, so that $\mathbf{x}= (d,-d)$.
From formula \eqref{eq:hain}, the Hain class is:
$$
\mathfrak{h}_1(d,-d) = \frac{d^2}{2} \left(\psi_1+\psi_2 \right).
$$ 

Because the Hain class coincides with the class of the double ramification cycle on compact type, the correction term must be supported on the irreducible divisor $\delta_0$:
\beq
DRC(d,-d) = \mathfrak{h}_1(d,-d) + \alpha \delta_0.
\eeq
One intersection theoretic computation is now sufficient to determine the unknown coefficient $\alpha$.
The double Hurwitz number is easily computed combinatorially as:
$$
H_1^2(d,-d)= \frac{1}{12} (d^3-d).
$$
On the other hand, the geometric formula  in \cite{cm:dhn} gives:

$$
H_1^2(d,-d)= r!(d\psi_1) DRC(d,-d) = d\psi_1(d^2(\psi_1+\psi_2)+ 2\alpha \delta_0).
$$
Performing the intersections one obtains $\alpha = -1/12$ and therefore
$$
DRC(d,-d) = \mathfrak{h}_1(d,-d) - \frac{1}{12}\delta_0 = \mathfrak{h}_1(d,-d) - \lambda_1.
$$

In genus higher than one, the only description - in terms of standard tautological classes - of a full compactification of the Hain class appears in work of Tarasca \cite{tarasca}, who computes the class of the admissible cover compactification in genus $2$ and for full ramification profiles over both  $0$ and $\infty$.

In \cite{GZ2}, Grushevski-Zakharov are able to extend the Hain class to the locus of curves with at most one non-separating node. They do so by computing the class of the zero section in a partial compactification of the moduli space of abelian varieties, and pulling it back via a naturally extended Abel-Jacobi map.
\begin{theorem}[\cite{GZ2}]
Let $\sM_{g,n}^{\leq1}$ be the partial compactification of the moduli space of curves parameterizing curves with at most one non-separating node. Define the coefficients:
$$
\eta_{a,b}:= \frac{(-1)^b (2b-1)!!}{2^{3b} a!}\sum_{x=0}^b \frac{2-2^{2x}B_{2x}}{(2b-2x-1)!!(2x-1)!! (b-x)!x!},
$$
where $B_{2x}$ denotes the Bernoulli number.
Then:
\beq
DRC(\mathbf{x})_{|\sM_{g,n}^{\leq1}} = \sum_{a+b=g} \eta_{a,b}\mathfrak{h}_g(\mathbf{x})^a\delta_0^b. 
\eeq
\end{theorem}

All the facts we stated about double Hurwitz numbers, the Hain class and its extension to the partial compactification  to $\sM_{g,n}^{\leq1}$ corroborate the following ``folklore" conjecture.

\begin{conjecture}\label{conj:drcpoly}
The class $DRC(\mathbf{x})$ is a (piecewise) even polynomial class in the $x_i$'s of degree $2g$, with coefficients in $R^g(\overline{\sM}_{g,n})$. 
\end{conjecture}

Further evidence for this conjecture is given in recent work of Buryak-Shadrin-Spitz-Zvonkine \cite{bssz:psi}: their definition of double ramification cycle  is the push-forward of the virtual class of the moduli space of relative stable maps. They compute generating functions for  intersection numbers of monomials of psi classes against the double ramification cycle and show that these intersection numbers behave compatibly with  the above Conjecture.

Denote by $\mathfrak{S}(z)$ the power series
$$
\mathfrak{S}(z) = \frac{\sinh(z/2)}
{z/2}
=
\sum_{k³0}\frac{z^{2k}}{2^{2k}(2k+1)!}
= 1 +
\frac{z^2}{24}
+
\frac{z^4}{1920}
+
\frac{z^6}{322560}
+ \ldots
$$

\begin{theorem}[\cite{bssz:psi}]
The following intersection formula holds: 
$$
\psi_s^{2g-3+n}
DRC_g(\mathbf{x}) = [z^
{2g}
]
\prod_{
i\not=s} \frac{\mathfrak{S}(x_iz)}{
\mathfrak{S}(z)}
,$$
where $[z^{2g}
]$ denotes the coefficient of $z^{2g}$ in the expansion of the analytic functions in the formula at $z=0$.
\end{theorem}
More complicated combinatorial formulas are produced for the intersection of the $DRC$ with an arbitrary monomial in psi classes of degree $2g-3+n$.

The most exciting recent progress in the quest for understanding $DRC(\mathbf{x})$ is a conjecture formulated by Pixton, following his work in collaboration with Pandharipande and Zvonkine on using CohFT techniques in order to produce tautological relations in the Chow ring of  the moduli space of curves. This work of Pixton is currently still unpublished; the following is a brief survey of the information communicated by Pixton in a mini-course  he gave on his work at the {\it Summer school in Gromov Witten theory} at Pingree Park, in the summer of 2014. 

Informally, a cohomological field theory (CohFT) is a collection of classes in the cohomology of the moduli space of curves  $\overline{\sM}_{g,n}$ parameterized by elements of the $n$-fold tensor product of a vector space $V$, called the Hilbert space of the theory. Such classes restrict  in a natural way to boundary divisors (see, for example \cite[\S 5]{pixcohft}, for the precise formulation of the axioms), providing a rich inductive structure that allows for reconstruction arguments. 
Pixton's initial observation is that the collection of  all Hain classes can be organized to be a ``cohomological field theory $\mathfrak{H}$ on $\sM_{g,n}^{ct}$". More precisely, define
$V = \bigoplus_{i\in \bZ}\langle e_i \rangle_\bC$, $\eta(e_i,e_j)= \delta_{i+j,0}$ and

$$
\mathfrak{H}_g(e_{x_1}\otimes \ldots \otimes e_{x_n})= 
\left\{
\begin{array}{cl}
\mathfrak{h}_g(\mathbf{x}) & \mbox{if $\sum x_i=0$}\\
0                                          & \mbox{else}
\end{array}
\right.
$$
It is simple to see that $\mathfrak{H}$ satisfies all CohFT axioms except for the splitting axiom along the irreducible divisor. Since we are restricting all classes to compact type - the failure of the splitting axiom along the irreducible divisor is irrelevant.

The next step is to expand the formula for the Hain class, and express it as a sum over trees, representing the dual graphs of the corresponding boundary strata; this appears (up to a factor of $2^g$) as the degree $g$ part of the following expression:
\beq \label{eq:pix1}
e^{2\Phi_{\mathbf{x}}^\ast ([\Theta])} = \sum_{T} \frac{1}{|\aut (T) |} \iota_\ast\left( 
\prod_{legs}   
e^{w_l^2\psi_l}
\prod_{edges}
 \frac{1-e^{w_e^2(\psi_t+\psi_e)}}{\psi_t+\psi_e} 
  \right)
 ,
\eeq
where  $T$ denotes  trees decorated like the appropriate tropical covers of the tropical line, the $w$ are the weights of legs and edges, and the $\psi_t$ and $\psi_e$ are psi classes at the tail and end of the corresponding edge.

Pixton proposes to extend the formula from the Hain class to the double ramification cycle by allowing the summation in \eqref{eq:pix1} to range over all suitably decorated graphs instead of over only trees. 

The natural generalization of \eqref{eq:pix1} over non-contractible graphs takes the following form:
\beq \label{eq:pix2}
 \sum_{\Gamma} \frac{1}{|\aut(\Gamma) | \dim(V)^{b_1(\Gamma)}} \iota_\ast\left( 
\prod_{legs}   
e^{w_l^2\psi_l}
\prod_{edges}
 \frac{1-e^{w_e^2(\psi_t+\psi_e)}}{\psi_t+\psi_e} 
  \right).
\eeq
Formula \eqref{eq:pix2} contains some divergent series and undefined terms. One problem arises from dividing by the dimension of the vector space $V= \bC\langle \bZ\rangle$, raised to the first betti number of the graph $\Gamma$. The other problem arises from the fact that when a graph $\Gamma$ contains  a loop at a vertex, then the balancing condition  at that vertex can be satisfied by adding an arbitrary pair of opposite integers on either end of the loop. Thus the sum over all  admissible decorations of $\Gamma$ ends up being a divergent series.
These issues are addressed as follows:
\begin{enumerate}
	\item The ring $\bZ$ is replaced by $\bZ/r\bZ$ so as to get rid of the divergencies: $\dim(V) = r$ and there are only finitely many decorations (mod $r$) arising from loops.
	\item The corresponding classes are shown to form a CohFT $\mathfrak{Pix}_g^r$.
	\item For $r$ sufficiently large, the classes are polynomial in $r$.
	\item The class $\mathfrak{Pix}_g(\mathbf{x})$ is defined by evaluating the polynomial class $\mathfrak{Pix}^r_g(\mathbf{x})$ at $r=0$. 
\end{enumerate}

Pixton is able to compute the class $\mathfrak{Pix}_g(\mathbf{x})$ and show that it satisfies the polynomiality conjecture, and agrees with the formula of Grushevski and Zakharov over the partial compactification of curves with at most one non-separating node. This brings substantial amount of support for the following exciting conjecture.

\begin{conjecture}[Pixton]
$$
DRC_g(\mathbf{x}) = 2^{-g}\mathfrak{Pix}_g(\mathbf{x}).
$$
\end{conjecture}
In breaking news, during the 2015 AMS algebraic geometry institute in Salt Lake City, a proof of Pixton's conjecture by Janda, Pandharipande, Pixton and Zvonkine was announced\cite{rahulslides}. The preprint is now available on the arxiv \cite{proofofpixton}.

\bibliographystyle{alpha}

\bibliography{biblio}
\end{document}